\documentclass[12pt] {article}
\usepackage{amsmath,amsthm}
\title{Non-commutative linear algebra and plurisubharmonic functions
of quaternionic variables. }
\date{}
\author{ Semyon Alesker
\\  { \normalsize Department of Mathematics, Tel Aviv University, Ramat Aviv}
 \\  { \normalsize 69978 Tel Aviv,
Israel }
\\ {\normalsize e-mail: semyon@post.tau.ac.il}}
\newcommand{\RR}{\mbox{\rm $~\vrule height6.5pt width0.5pt
depth0.3pt\!\!$R}}

\newcommand{\CC}{\mbox{\rm $~\vrule height6.5pt width0.5pt
depth0.3pt\!\!$C}}

\newcommand{\HH}{\mbox{\rm $~\vrule height6.5pt width0.5pt
depth0.3pt\!\!$H}}

\def\eps{\varepsilon}

\def\lam{\lambda}

\def\str{\longrightarrow}

\def\d{det}
\def\qed { Q.E.D. }

\def\rp{\RR_{>0}}
\def\rnn{ \RR _{\geq 0}}
\def\hab{ \HH^*_{ab}}
\def\mii{M'_{II}(A)}
\def\miin{M'_{II}(A_1)}
\def\mij{M'_{IJ}(A)}

\def\mjj{M'_{JJ}(A)}

\def\mpij{M'_{I\cap J,I\cap J}(A)}

\def\psh{plurisubharmonic }

\def\db{\frac{\partial}{\partial \bar q}}
\def\dq{\frac{\partial}{\partial  q}}

\def\dfq{\frac{\partial ^2 f}{\partial \bar q_i \partial q_j}}
\def\ih{\int_{\HH ^n}}
\swapnumbers
\newtheorem{theorem}{Theorem}[subsection]
\newtheorem{corollary}[theorem]{Corollary}
\newtheorem{lemma}[theorem]{Lemma}
\newtheorem{proposition}[theorem]{Proposition}
\newtheorem{claim}[theorem]{Claim}
\theoremstyle{definition}
\newtheorem{example}[theorem]{Example}
\newtheorem{definition}[theorem]{Definition}
\newtheorem{remark}[theorem]{Remark}
\begin{document}
\maketitle
\begin{abstract}
We remind known and establish new properties of the Dieudonn\'e
and Moore determinants of quaternionic matrices. Using these linear
algebraic results we develop a basic theory of plurisubharmonic
functions of quaternionic variables.
\end{abstract}
\setcounter{section}{-1}
\section{Introduction.}
The main point of this paper is that
in quaternionic algebra and analysis there exist
structures which have analogues over
the fields of real and complex
numbers, but should reflect different phenomena.

The algebraic part is discussed in Section 1. There we remind the
notions of the Moore and Dieudonn\'e determinants of quaternionic
matrices. It turns out that (under appropriate normalization) the
Dieudonn\'e determinant behaves exactly like the absolute value of
the usual determinant of real or complex matrices from all points
of view (algebraic and analytic). Let us state some of its
properties discussed in more details in Subsection 1.2. Let us
denote by $M_n(\HH)$ the set of all quaternionic $n\times n$-
matrices. The Dieudonn\'e determinant $D$ is defined on this set
and takes values in non-negative real numbers:
$$D:M_n(\HH) \str \rnn$$
(see Definition 1.2.2). Then one has the following (known) results
(see  Theorems 1.2.3 and 1.2.4 below and references given
at the beginning of Section 1):

{\bf Theorem.} {\itshape

(i) For any complex $n \times n$- matrix $X$
considered as quaternionic matrix the Dieudonn\'e
determinant $D(X)$ is equal to the absolute value
of the usual determinant of $X$.

(ii) For any quaternionic matrix $X$
$$D(X^*)=D(X),$$
where  $X^*$ denotes
the quaternionic conjugate matrix of $X$.

(iii) $D(X\cdot Y) =D(X) D(Y).$
}

{\bf Theorem 1.2.5.} {\itshape
Let $A= \left[ \begin{array}{ccc}
a_{11}&\dots&a_{1n}\\
\multicolumn{3}{c}{\dotfill} \\
a_{n1}&\dots&a_{nn}
\end{array}
\right] $ be a quaternionic matrix.
Then $$D(A) \leq \sum_{i=1}^n |a_{1i}| D(M_{1i}).$$
Similar inequalities hold for any other row or column.
}

(In this theorem $|a|$ denotes the absolute value of a quaternion
$a$, and $M_{pq}$ denotes the minor of the matrix $A$ obtained
from it by deleting the $p$-th row and $q$-th column).

In a sense, the Dieudonn\'e determinant provides the theory of
{\itshape absolute value } of determinant. However it is not
always sufficient and we loose most of the algebraic properties of
the usual determinant. The notion of Moore determinant provides
such a theory, but only on the class of quaternionic {\itshape
hyperhermitian } matrices. Remind that a square quaternionic
matrix $A=(a_{ij})$ is called hyperhermitian if its quaternionic
conjugate $A^*=A$, or explicitly $a_{ij}=\overline {a_{ji}}$.  The
Moore determinant
 denoted by $det$
is defined on the class of all hyperhermitian matrices
and takes real values. (The Moore determinant is defined
in Subsection 1.1 after Theorem 1.1.8).
 The important advantage of it with respect to the Dieudonn\'e
determinant is that it depends polynomially on the entries of a
matrix; it has already all the
 algebraic and analytic properties of the usual determinant of real symmetric
and complex hermitian matrices. Let us state some of them referring for the
details to Subsection 1.1
(again, the references are given at the beginning of Section 1).

{\bf Theorem 1.1.9.} {\itshape

(i) The Moore determinant of any
complex hermitian matrix considered as quaternionic hyperhermitian matrix
is equal to its usual determinant.

(ii) For any hyperhermitian matrix $A$ and any quaternionic matrix $C$
$$det (C^*AC)= detA \cdot det(C^*C).$$
}

{\bf Examples.}

(a) Let $A =diag(\lam_1, \dots, \lam _n)$ be a diagonal matrix
with real $\lam _i$'s. Then $A$ is hyperhermitian
and the Moore determinant
$detA= \prod _i \lam_i$.

(b)  A general hyperhermitian $2 \times 2$ matrix $A$ has the form
 $$ A=  \left[ \begin {array}{cc}
                     a&q\\
                \bar q&b\\
                \end{array} \right] ,$$
where $a,b \in \RR, \, q \in \HH$. Then its Moore determinant
is equal to
$det A =ab - q \bar q$.

Next, in terms of the Moore determinant one can prove the
generalization of the classical Sylvester criterion of positive
definiteness of hyperhermitian matrices (Theorem 1.1.13). In terms
of the Moore determinant one can introduce the notion of the mixed
discriminant and to prove the analogues of Aleksandrov's
inequalities for mixed discriminants (Theorem 1.1.15 and Corollary
1.1.16).

The (well known) relation between the Dieudonn\'e and Moore determinants
is as follows: for any hyperhermitian matrix $X$
$$D(X)=|det X|.$$
In Section 1 we prove some additional properties of the
Dieudonn\'e and Moore determinants; they are used in Section 2.

Note that the Dieudonn\'e determinant was introduced originally by J. Dieudonn\'e
in  \cite{dieudonne} (see also \cite{artin} for his theory).
 It can be defined for arbitrary (non-commutative) field.
On more modern language this result can be formulated as a
computation of the $K_1$- group of a non-commutative field (see
e.g. \cite{rosenberg}). Note also that there is a more recent
theory of non-commutative determinants (or quasideterminants) due
to I. Gelfand and V. Retakh generalizing in certain direction the
theory of the Dieudonn\'e determinant. First it was introduced in
\cite{gelfand-retakh1}, see also \cite{gelfand-retakh2},
\cite{gelfand-retakh3}, \cite{gelfand-gelfand-retakh-wilson} and
references therein for further developments and applications. In
the recent preprint \cite{gelfand-retakh-wilson} Gelfand, Retakh,
and Wilson have discovered that the formulas for quasideterminants
of quaternionic matrices can be significantly simplified. They
also understood the relation between the theory of
quasideterminants and the Moore determinant. We would also like to
mention a different direction of a development of the {\itshape
quaternionic} linear algebra started by D. Joyce \cite{joyce} and
applied by himself to hypercomplex algebraic geometry. We refer
also to D. Quillen's paper \cite{quillen} for further
investigations in that direction. Another attempt to understand
the quaternionic linear algebra from the topological point of view
was done in \cite{alesker2}.

Section 2 of this paper develops the basic theory of plurisubharmonic
functions of quaternionic variables on $\HH ^n$. It uses in essential
way the linear algebraic
results of Section 1. This theory is parallel to the classical theories
of convex functions on $\RR^n$  and plurisubharmonic functions on $\CC ^n$.

The formal definition is as follows (for more discussion see
Subsection 2.1).

{\bf Definition.} {\itshape A  real valued function
$$f:\HH ^n \str \RR$$
 is called quaternionic \psh
if it is upper semi-continuous and
 its restriction to any right quaternionic line is subharmonic
(in the usual sense).}
\newline
We refer to Subsection 2.1 where we remind the relevant notions.

 In this form this definition
 was suggested by G. Henkin \cite{henkin}. For the class of continuous \psh functions
this definition is different but equivalent (by Proposition 2.1.6 below) to the original author's definition.

{\bf Remark.} On $\HH^1$ the class of plurisubharmonic functions
coincides with the class of subharmonic functions. In this case
all the results of this paper are reduced to the classical
properties of subharmonic functions in $\RR^4$.

Let us describe the main results on \psh functions we prove.
We will write a quaternion $q$ in the usual form
$$q= t+ x\cdot i +y\cdot j+ z\cdot k ,$$
where $t,\, x,\, y,\, z$ are real numbers, and $i,\, j,\, k$
satisfy the usual relations
$$i^2=j^2=k^2=-1, \, ij=-ji=k,\, jk=-kj=i, \, ki=-ik=j.$$
Let us introduce the differential operators $\db$ and $\dq$ as follows:
$$\db f:=\frac{\partial f}{\partial  t}  +
i \frac{\partial f}{\partial  x} +
j \frac{\partial f}{\partial  y} +
k \frac{\partial f}{\partial  z}  , \mbox { and }$$
$$\dq f:=\overline{ \db \bar f}=
\frac{\partial f}{\partial  t}  -
 \frac{\partial f}{\partial x}  i-
 \frac{\partial f}{\partial  y} j-
\frac{\partial f}{\partial  z}  k.$$

{\bf Remarks.} (a) The operator $\db$ is called sometimes the
Cauchy-Riemann-Moisil-Fueter operator since it was introduced by
Moisil in \cite{moisil} and used by Fueter \cite{fueter1},
\cite{fueter2} to define the notion of quaternionic analyticity.
For further results on quaternionic analyticity we refer e.g. to
\cite{brackx-delanghe-sommen}, \cite{palamodov},
 \cite{pertici}, \cite{sudbery}, and for applications to mathematical
physics to \cite{gursey-tze}.
   Another used name for this
 operator is Dirac-Weyl operator. But in fact it was used earlier
 by J.C. Maxwell  in \cite{maxwell}, vol. II, pp.570-576, where he
 has applied the quaternions to electromagnetism.

 (b) Note that $$\db= \frac{\partial}{\partial t}+ \nabla,$$
 where $\nabla= i\frac{\partial}{\partial x}+j\frac{\partial}{\partial
 y}+ k\frac{\partial}{\partial z}$. The operator $\nabla$ was
 first introduced by W.R. Hamilton in \cite{hamilton}.

 (c) In quaternionic analysis one considers a right version of the
 operators $\db$ and $ \dq$ which are denoted respectively by
 $\overset{\leftarrow}{\db}$ and $ \overset{\leftarrow}{\dq}$. The
 operator $ \overset{\leftarrow}{\dq}$ is related to
 $\overset{\leftarrow}{\db}$ by the same formula as $\dq$ is
 related to $\db$, and $\overset{\leftarrow}{\db}$ is defined as
 $$\overset{\leftarrow}{\db}f :=\frac{\partial f}{\partial  t}  +
 \frac{\partial f}{\partial  x}i +  \frac{\partial f}{\partial y}j
+  \frac{\partial f}{\partial  z}k .$$ For a real valued function
$f$ the derivatives $\frac{\partial ^2 f}{\partial q_i \partial
\bar q_j}$ and $\frac{\partial ^2 f}{\overset{\leftarrow}{\partial
q_i} \overset{\leftarrow}{\partial \bar q_j}}$ are quaternionic
conjugate to each other.

First one has a simple

{\bf Proposition 2.1.6.}{\itshape
A real valued twice continuously
differentiable function $f$ on the domain $\Omega \subset \HH ^n$
is quaternionic
plurisubharmonic if and only if at every point $q \in \Omega$
the matrix $(\frac{\partial ^2 f}{\partial \bar q_i \partial  q_j})(q)$
is non-negative definite.
}

Note that the matrix in the statement of proposition is
quaternionic hyperhermitian. The more important thing is that
in analogy to the real and complex cases one can define
for any continuous quaternionic \psh function $f$ a
non-negative measure $det(\frac{\partial ^2 f}{\partial \bar q_i \partial q_j})(q)$,
where $det$ denotes the Moore determinant (this measure is obviously
defined for smooth $f$). We prove the following continuity result.

{\bf Theorem 2.1.11.} {\itshape
Let $\{f_N\}$ be sequence of continuous
quaternionic \psh function in a domain $\Omega \subset \HH^n$.
Assume that this sequence converges uniformly
on compact subsets
to a function $f$. Then $f$ is continuous quaternionic \psh
function. Moreover the sequence of measures
$det(\frac{\partial ^2 f_N}{\partial \bar q_i \partial  q_j})$
weakly converges to the measure $det(\dfq)$.
}

The proofs of analogous results in real and complex cases can be
found in \cite{aubin}, where the exposition of this topic follows
the approach of Chern-Levine-Nirenberg \cite{chern-levine-nirenberg} and
  Rauch-Taylor \cite{rauch-taylor}. For the complex
case we refer to the classical book by P. Lelong \cite{lelong}. In
generalizations of these results to the quaternionic situation the
large part of the difficulties comes from linear algebra since the
technique of working with the Moore determinant is not
sufficiently developed. For instance there is no formula of
decomposition of the Moore determinant in row or column, and thus
one should use some more tricky manipulations.

Next we would like to state a result on existence and uniqueness
of solution of the Dirichlet problem
 for quaternionic Monge-Amp\`ere equation (to be defined).
 In this paper we prove only the uniqueness part; the existence
 is proved in author's paper \cite{alesker3}.

{\bf Definition.}{\itshape  An open bounded domain $\Omega\subset
\HH^n$ with a smooth boundary $\partial \Omega$ is called strictly
pseudoconvex if for every point $z_0\in \partial \Omega$ there
exists a neighborhood ${\cal O}$ and a smooth strictly \psh
function $h$ on ${\cal O}$ such that $\Omega \cap {\cal O}= \{ h<0
\}$ and $\nabla h(z_0) \ne 0$.}
\newline

The next result is quaternionic analogue of the results on
Dirichlet problem for real and complex Monge-Amp\`ere equations.
The real case was solved by Aleksandrov \cite{aleksandrov2}, and
the complex one by Bedford and Taylor  \cite{bedford-taylor}.

{\bf Theorem.} {\itshape Let $\Omega$ be a strictly pseudoconvex
bounded domain in $\HH^n$. Let $\phi$ be a continuous real valued
function on the boundary $\partial \Omega$. Let $f$ be a
continuous function on the closure $\bar \Omega$, $f \geq 0$. Then
there exists a unique continuous on $\bar \Omega$ \psh function
$u$ such that
$$det(\frac{\partial ^2 u}{\partial \bar q_i \partial  q_j})
=f \mbox{ and}$$
$$u=\phi \mbox{ on } \partial \Omega.$$
}

The uniqueness part in this theorem is an immediate consequence of
the following {\itshape minimum principle} which is proved in
Subsection 2.2.

{\bf Theorem 2.2.1.} {\itshape Let $\Omega$ be a bounded  open set
in $\HH^n$. Let $u,\, v$ be continuous functions on $\bar \Omega$
which are \psh in $\Omega$. Assume that
$$det(\frac{\partial ^2 u}{\partial \bar q_i \partial q_j})\leq
det(\frac{\partial ^2 v}{\partial \bar q_i \partial  q_j}) \mbox{
in } \Omega.$$ Then
$$min\{u(z)-v(z)|z\in \bar \Omega\} =min\{u(z)-v(z)|z\in \partial \Omega\}.$$ }

The proof of this theorem closely follows the argument of Bedford
and Taylor \cite{bedford-taylor} (Theorem A).

In appendix to this paper we prove the injectivity of Radon
transform over quaternionic subspaces in the affine space
$\HH ^n$. Probably this result is not new. It is included here
since it was used in the proof of Lemma 2.1.7, and we could not find
a reference.

{\bf Acknowledgements.} We express our gratitude to J. Bernstein,
P. Biran, G. Henkin, D. Kazhdan, V. Palamodov, L. Polterovich, and
M. Sodin
 for useful discussions. We thank also E. Shustin
who has shown us the reference \cite{aslaksen}.

\section { Linear algebra.}
In this section we remind the construction and basic properties of
the Dieudonn\'e and Moore determinants and investigate further
their properties. Part of them will be used in the next sections
of this paper. For a survey of quaternionic determinants and
references we refer to \cite{aslaksen}.

First of all remind that over any noncommutative field
there exist usual notions of vector spaces over the field
(however one should distinguish between left and right ones),
their dimension, basis etc. (see e.g. \cite{artin}). However there is
no construction of quaternionic determinant which would have all
the properties of the determinant over commutative field.
We are going to discuss this problem in this section.
We will discuss only right vector spaces. The case of left ones
can be considered similarly. Many results of Section 1 are a folklore.
Theorems 1.1.8, 1.1.9, 1.1.4 are not new. We refer for the proofs to
\cite{barnard-moore}, \cite{chen1},\cite{chen2},
\cite{dyson1},\cite{dyson2}, \cite{jacobson1}, \cite{jacobson2},
\cite{mehta1}, \cite{mehta2}, \cite{moore}, \cite{piccinni},
\cite{vanpraag1}, \cite{vanpraag2}.

\subsection {Hyperhermitian forms and the Moore determinant.}
Let $V$ be a right vector space over quaternions.

\begin{definition}
  A {\itshape hyperhermitian semilinear form}
on $V$ is a map $ a:V \times V \str \HH$ satisfying the following properties:

(a) $a$ is additive with respect to each argument;

(b) $a(x,y \cdot q)= a(x, y) \cdot q$ for any $x,y \in V$ and any $q\in \HH$;

(c) $a(x,y)= \overline{a(y,x)}$.
\end{definition}

\begin{example} Let $V= \HH ^n$ be the standard coordinate space
considered as right vector space over $\HH$. Fix a {\itshape
hyperhermitian} $n \times n$-matrix $(a_{ij})_{i,j=1}^{n}$, i.e.
$a_{ij} =\bar a_{ji}$, where $\bar x$ denotes the usual
quaternionic conjugation of $x$. For $x=(x_1, \dots, x_n), \,
y=(y_1, \dots, y_n)$ define
$$A(x,y) = \sum _{i,j} \bar x_i a_{ij} y_j$$
(note the order of the terms!).
Then $A$ defines hyperhermitian semilinear form on $V$.
\end{example}
In general one has the following standard claims.

\begin{claim} Fix a basis in a finite dimensional right quaternionic
vector space $V$. Then there is a natural bijection between
hyperhermitian semilinear forms on $V$ and $n \times
n$-hyperhermitian matrices.
\end{claim}
This bijection is in fact described in previous Example 1.1.2.

\begin{claim} Let $A$ be the matrix
of the given hyperhermitian form in the given basis. Let
$C$ be transition matrix from this basis to another one.
Then the matrix $A'$ of the given form in the new basis
is equal $$A' =C^* AC .$$
\end{claim}
\begin{remark}
Note that for any hyperhermitian matrix $A$ and for any matrix $C$
the matrix $C^* AC$ is also hyperhermitian. In particular the
matrix $C^* C$ is always hyperhermitian.
\end{remark}

\begin{definition} A hyperhermitian semilinear form $a$
is called {\itshape positive definite} if $a(x,x)>0$ for any non-zero vector
$x$.
 \end{definition}

Let us fix on our quaternionic right vector space $V$ a
positive definite hyperhermitian form $( \cdot , \cdot )$.
The space with  fixed such a form will be called
{\itshape hyperhermitian} space.

For any quaternionic linear operator $\phi: V\str V$ in
hyperhermitian space
one can define the adjoint operator $\phi ^* :V \str V$
in the usual way, i.e.
$(\phi x,y)= (x, \phi ^* y)$ for any
$x,y \in V$. Then if one fixes an orthonormal basis
in the space $V$ then the operator $\phi$ is selfadjoint
if and only if its matrix in this basis is hyperhermitian.

\begin{claim}
For any  selfadjoint operator in a hyperhermitian space
there exists an orthonormal basis such that
its matrix in this basis is diagonal and real.
\end{claim}
The proof is standard. Now we are going to define the Moore determinant
of  hyperhermitian matrices. The definition below
is different from the original one \cite{moore} but equivalent to it.

First note that every hyperhermitian $n \times n$- matrix $A$
defines a hyperhermitian semilinear form on the coordinate space
$\HH ^n$. It also can be considered as a {\itshape symmetric}
bilinear form on $\RR ^{4n}$ (which is the realization of $\HH
^n$). Let us denote its $4n \times 4n$- matrix by ${}^{\RR} A$.
Let us consider the entries of $A$ as formal variables (each
quaternionic entry corresponds to four commuting real variables).
Then $det ({}^{\RR} A)$  is a homogeneous polynomial of degree
$4n$ in $n(2n-1)$ real variables. Let us denote by $Id$ the
identity matrix.
 One has the following result.

\begin{theorem}
There exists a polynomial $P$ defined on the space of all
hyperhermitian $n \times n$-matrices such that for any
hyperhermitian $n \times n$-matrix $A$ one has $det({}^{\RR} A)=
P^4(A)$ and $P(Id)=1$. $P$ is defined uniquely by these two
properties. Furthermore $P$ is homogeneous of degree $n$ and has
integer coefficients.
\end{theorem}
Thus for any hyperhermitian matrix $A$ the value $P(A)$ is a real
number, and it is called the {\itshape Moore determinant} of the
matrix $A$. The explicit formula for the Moore determinant  was
given by Moore \cite{moore} (see also \cite{aslaksen}). From now
on the Moore determinant of a matrix $A$ will be denoted by $det
A$. This notation should not cause any confusion with the usual
determinant of real or complex matrices due to part (i) of the
next theorem.

\begin{theorem}

(i) The Moore determinant of any
complex hermitian matrix considered as quaternionic hyperhermitian matrix
is equal to its usual determinant.

(ii) For any hyperhermitian matrix $A$ and any matrix $C$
$$det (C^*AC)= detA \cdot det(C^*C).$$
\end{theorem}
\begin{example}

(a) Let $A =diag(\lam_1, \dots, \lam _n)$ be a diagonal matrix
with real $\lam _i$'s. Then $A$ is hyperhermitian
and the Moore determinant
$detA= \prod _i \lam_i$.

(b)  A general hyperhermitian $2 \times 2$ matrix $A$ has the form
 $$ A=  \left[ \begin {array}{cc}
                     a&q\\
                \bar q&b\\
                \end{array} \right] ,$$
where $a,b \in \RR, \, q \in \HH$. Then
$det A =ab - q \bar q$.
\end{example}
Let us introduce more notation. Let $A$ be any hyperhermitian $n
\times n$- matrix. For any non-empty subset $I \subset \{1, \dots,
n\}$  the minor $M_I(A)$ of $A$ which is obtained by deleting the
rows and columns with indexes from the set $I$, is clearly
hyperhermitian.
For $I=\{1, \dots , n \}$ let $det M_{\{1, \dots,n \}} =1$.
\begin{proposition}
For any hyperhermitian $n \times n$-matrix $A$ and any diagonal
real matrix $T= \left[ \begin{array}{ccc}
                     t_1&        &    0 \\
                        & \ddots & \\
                      0 &        &t_n \\
            \end{array} \right]$
$$det (A+T) =
\sum _{I\subset \{1, \dots ,n\} } (\prod _{i\in I} t_i) \cdot det M_I(A).$$
In particular
$$det(A+t \cdot Id) =\sum _{I\subset \{1, \dots ,n\} } t^{|I|} \cdot det M_I(A),$$
where $|I|$ denotes the cardinality of the set $I$.
\end{proposition}

{\bf Remark.} Clearly this formula is true for {\itshape
arbitrary} $n\times n$-matrix $A$ over a {\itshape commutative}
field.

{\bf Proof.} Fix a hyperhermitian matrix $A$. It is clear that
$det (A+T)$ is a polynomial in $t_1, \dots, t_n$ of degree $n$.
Since $$A+ \left[ \begin{array}{cccc}
               t_1&   &   &0\\
                  &t_2&   & \\
                  &   &\ddots & \\
                0 &   &   &t_n\\
                  \end{array} \right] = A+
\left[ \begin{array}{cccc}
                 0&   &   &0\\
                  &t_2&   & \\
                  &   &\ddots & \\
                0 &   &   &t_n\\
                  \end{array} \right]+
\left[ \begin{array}{cccc}
               t_1&   &   &0\\
                  &0  &   & \\
                  &   &\ddots & \\
                0 &   &   &0\\
                  \end{array} \right]$$
one can apply induction in $n$  to show that if
$A= \left[ \begin{array}{c|c}
                 a_{11} &\ast \\ \hline
                  \ast  &B\\
                   \end{array} \right] ,$
where $a_{11} \in \RR$, and $B$ is a hyperhermitian
$(n-1)\times (n-1)$- matrix.
Set $$f(t):= det \left( A+  \left[ \begin{array}{cccc}
                 t&   &   &0\\
                  &0  &   & \\
                  &   &\ddots & \\
                0 &   &   &0\\
                  \end{array} \right] \right).$$
It is sufficient to show that
$ f(t)= det A + t\cdot detB $.
Clearly $f(0)= det A$. Let $k$ denote the degree of the polynomial $f$.
Using Theorem 1.1.9(ii) one gets
$f(t)=$
$$ t^k \cdot det \left(
    \left[ \begin{array}{cccc}
          t^{-k/2}& & &0 \\
                  &1& & \\
                  & &\ddots& \\
                0 & & & 1
             \end{array} \right]
  \left( A+  \left[ \begin{array}{cccc}
                 t&   &   &0\\
                  &0  &   & \\
                  &   &\ddots & \\
                0 &   &   &0\\
                  \end{array} \right] \right)
    \left[ \begin{array}{cccc}
          t^{-k/2}& & & 0\\
                  &1& & \\
                  & &\ddots& \\
                0 & & & 1
             \end{array} \right]
\right)= $$
$$t^k det \left[ \begin{array}{c|ccc}
a_{11} t^{-k} +t^{-k+1}&
a_{12} t^{-k/2}& \dots & a_{1n}t^{-k/2} \\ \hline
                            \begin{array}{c}
                            a_{21}t^{-k/2} \\
                                            \vdots \\
                            a_{n1}t^{-k/2} \end{array} & & B&
\end{array}
\right].$$
If $k>1$ then
$$\frac{f(t)}{t^k} \str det \left[ \begin{array}{c|c}
0&0\\ \hline
0&B
\end{array}
\right]=0$$
when $t \str \infty$.
Hence $k = 1$ and
$$\frac{f(t)}{t} \str det \left[ \begin{array}{c|c}
1&0\\ \hline
0&B
\end{array}
\right]= det B .$$
\qed
\begin{lemma}
Let $A$ be  a non-negative (resp. positive) definite
hyperhermitian matrix.
Then $det A \geq 0 \, (\mbox{ resp. } det A >0)$.
\end{lemma}
{\bf Proof.} Let us prove it under the assumption that $A$
is positive definite.
 By Claim 1.1.7 there exists a matrix $C \in Sp(n)$ (i.e. $C^* C =Id$)
such that $A=C^* \left[ \begin{array}{ccc}
                        \lam _1 &   &0 \\
                                &\ddots &\\
                             0  &   &\lam _n
                                \end{array} \right] C $
with
$ \lam _i \in \RR$. Since $A$ is positive definite, $\lam _i >0$ for all $i$.
By Theorem 1.1.9(ii) $det A= det (C^*C) \prod \lam _i =\prod \lam _i >0$.
\qed

The following theorem is a quaternionic generalization of the standard
Sylvester criterion.
\begin{theorem}[Sylvester criterion]
A hyperhermitian $n \times n$- matrix $A$ is positive definite if
and only if $M_{\{ i+1, \dots , n\}}(A) >0$ for any $0\leq i \leq
n$.
\end{theorem}
{\bf Proof.} The necessity follows from Lemma 1.1.12. Let
us prove sufficiency by induction in $n$. For $n=1$ the statement
is trivial. Assume $n>1$. Let $$ A= \left[ \begin{array}{c|ccc}
a_{11}&a_{12}&\dots&a_{1n} \\ \hline
a_{21}& & & \\
\vdots& & B& \\
a_{n1}& & &
\end{array}
\right].$$
Consider the matrix
$$ U= \left[ \begin{array}{ccccc}
1& -a_{12}/a_{11} & -a_{13}/a_{11} & \dots &-a_{1n}/a_{11} \\
0&1&0& \dots &0 \\
\multicolumn{5}{c}{\dotfill} \\
0& 0& 0& \dots &1 \\
\end{array}
\right] .$$
Then the matrix $A':= U^* AU$ has the form
$$A'= \left[ \begin{array}{cccc}
a_{11}&0&\dots&0 \\
0& & & \\
\vdots& & B'& \\
0& & &
\end{array}
\right],$$ where $B'$ is a hyperhermitian matrix. Moreover for any
$1 \leq i \leq n$ one has
$$det M_{\{i+1, \dots ,n \}} (A') =
det M_{\{i+1, \dots ,n\}} (A).$$
Indeed let us check it for $i=n$
(for $i<n$ the proof will be the same since the matrix $U$ is triangular).
Namely let us show that $det A'= det A$.
By Theorem 1.1.9 (ii)
 $det A'= det A \cdot det(U^*U)$. However using Theorem 1.1.8
and unipotence of $U$
it is  easy to see that $det(U^*U)=1$.
Hence the matrix $B'$ is positive definite by the induction
assumption. Then $A'$ is positive definite, and hence $A$ is as well.
\qed

Let us define now the mixed discriminant of hyperhermitian
matrices in analogy with the case of real symmetric matrices \cite{alex}.
\begin{definition}
Let $A_1, \dots ,A_n$ be hyperhermitian $n \times n$- matrices.
Consider the homogeneous polynomial in real variables $\lam _1
,\dots , \lam _n$ of degree $n$ equal to $det(\lam_1 A_1 + \dots +
\lam_n A_n)$. The coefficient of the monomial $\lam_1 \cdot \dots
\cdot \lam_n$ divided by $n!$ is called the {\itshape mixed
discriminant} of the matrices  $A_1, \dots ,A_n$, and it is
denoted by $\d(A_1, \dots ,A_n)$.
\end{definition}
Note that the mixed discriminant is symmetric with respect to all variables,
and linear with respect to each of them, i.e.
$$\d (\lam A_1' +\mu A_1'', A_2, \dots, A_n )=
\lam \cdot \d( A_1', A_2, \dots, A_n ) + \mu \cdot \d(A_1'', A_2,
\dots, A_n )$$ for any {\itshape real} $\lam , \, \mu$. Note also
that $\d(A, \dots, A)=det A$. We will prove the following
generalization of Aleksandrov's inequalities for mixed
discriminants \cite{alex} (though the proof will be very close to
the original one).
\begin{theorem}
(i) The mixed discriminant of positive (resp. non-negative) definite matrices
 is positive (resp. non-negative).

(ii) Fix positive definite hyperhermitian $n \times n$- matrices
$A_1, \dots, A_{n-2}$. On the real linear space
of hyperhermitian $n \times n$- matrices consider
the bilinear form $$B(X,Y):= \d(X,Y, A_1, \dots, A_{n-2}).$$
Then $B$ is non-degenerate quadratic form, and
its signature has one plus and the rest are minuses.
\end{theorem}
\begin{corollary}
 Let
$A_1, \dots, A_{n-1}$ be positive definite hyperhermitian $n \times n$- matrices.
Then for any hyperhermitian matrix $X$
\begin{equation}
\d(A_1, \dots, A_{n-1}, X)^2 \geq \d(A_1, \dots, A_{n-1},A_{n-1})
\cdot \d(A_1, \dots , A_{n-2}, X,X) ,
\end{equation}
and the equality is satisfied if and only if the matrix $X$ is proportional
to $A_{n-1}$.
\end{corollary}
{\bf Proof} of Corollary 1.1.16. By Theorem 1.1.15 (i) we get
$$\d(A_1, \dots, A_{n-1},A_{n-1}) >0.$$ Let
$$\lam =\frac{\d(A_1, \dots, A_{n-1},X)}{\d(A_1, \dots, A_{n-1},A_{n-1})}.$$
Let $X'=X- \lam A_{n-1}$. Then clearly $\d(A_{1}, \dots ,A_{n-1},
X')=0$. In the notation of Theorem 1.1.15 it means that
$B(A_{n-1}, A_{n-1}) >0$ and $B(A_{n-1}, X') =0$. But the form $B$
has just one plus. Hence $B(X',X') \leq 0$, and the equality is
satisfied if and only if $X'=0$. Developing $B(X',X')$ one gets
inequality (1). The equality case follows as well. \qed

{\bf Proof} of Theorem 1.1.15. (1) Let us prove the first part
using induction in $n$. The case $n=1$ is trivial.
Assume that $n>1$. Let $A_1, \dots ,A_n$ be positive definite
hyperhermitian matrices.
By Claim 1.1.7 and Theorem 1.1.9 (ii)
 we can assume that the matrix $A_n$ is diagonal, i.e.
$A_n= \left[ \begin{array}{ccc}
t_1& &0 \\
 &\ddots& \\
0& &t_n \\
\end{array}
\right] $, and $t_i$'s are positive. By Proposition 1.1.11
$$det (\lam _1 A_1 + \dots +\lam_{n-1} A_{n-1} +
\lam_n A_n)=$$
$$ \sum_{I \subset \{1, \dots, n\} }
(\prod _{i \in I} \lam_n t_i) \cdot det M_I (\lam_1 A_1 + \dots
+\lam_{n-1} A_{n-1}) .$$
Since all the diagonal minors of positive definite matrix
are positive definite and since $t_i >0$ the
assumption of induction implies the statement.

(2) Let us prove the second part of the theorem,
i.e. that $B$ is non-degenerate.
First let us prove it for $n=2$. Assume $X_0$ belongs to the kernel
of $B$, i.e. $B(X,X_0)=0$ for every $X$. One can assume that $X_0$
is diagonal: $X_0= \left[ \begin{array}{cc}
t_1&0\\
0&t_2
\end{array}
\right]$. For any $X=\left[ \begin{array}{cc}
x_1&0\\
0&x_2
\end{array}
\right]$ with real $x_1, \, x_2$ one has
 $2 \d(X,X_0)= t_1 x_2 +t_2 x_1 =0$. Hence $t_1=t_2=0$.
Thus the form $B$ is non-degenerate. Now, clearly
$B(Id, Id)=1 >0$. Assume that $X\ne 0$ is orthogonal to $Id$ with respect to $B$,
i.e. $B(X,Id)=0$. It remains to show that $B(X,X) <0$.
By Claim 1.1.7 we can assume that $X$ is diagonal,
$X=\left[ \begin{array}{cc}
x_1&0\\
0&x_2
\end{array}
\right]$. Then $2 B(X,Id)=x_1+x_2 =0$. But
$ B(X,X)=x_1 x_2=-x_1^2 <0$.

Let us assume that $n>2$. Assume also that the theorem is true for
matrices of size at most $n-1$. Let us prove first that the form
$B$ is non-degenerate. Assume that $X_0$ belongs to the kernel of
$B$. Since $A_{n-2}$ is positive definite, by Claim 1.1.7 one can
assume that the matrix $A_{n-2}$ is equal to $Id$ and $X_0$ is
diagonal. For $1\leq i \leq n$ and for $n\times n$- matrices $C_1,
\dots , C_{n-1}$ let us denote by $\d(C_1, \dots , C_{n-1})_i$ the
mixed discriminant of $(n-1)\times (n-1)$-matrices obtained from
$C_j$'s by deleting the $i$-th row and the $i$-th column. Let $T=
\left[
\begin{array}{ccc}
t_1& &0\\
 &\ddots& \\
0& &t_n
\end{array}
\right] $.
Using Proposition 1.1.11 one can easily see that
\begin{equation}
0=c \cdot \d(A_1, \dots ,A_{n-2}, T,X_0)= \sum _{i=1}^n
t_i \d(A_1, \dots ,A_{n-2}, X_0)_i ,
\end{equation}
where $c>0$ is a normalizing constant.
Hence
$\d(A_1, \dots ,A_{n-2}, X_0)_i =0$ for all $i$.
By the induction assumption and Corollary 1.1.16
(which is also satisfied for matrices of size $n-1$)
\begin{equation}
\d(A_1, \dots ,A_{n-3}, X_0,X_0)_i \leq 0 \mbox{ for } i=1, \dots ,n
\end{equation}
with equalities if and only if the matrix $X_0$  vanishes.
Since $A_{n-2}=Id$ and $X_0$ belongs to
the kernel of $B$ the equality analogous to (2) implies that
 $$0=c \cdot \d(A_1, \dots ,A_{n-3}, A_{n-2},X_0,X_0)= \sum _{i=1}^n
 \d(A_1, \dots ,A_{n-3}, X_0,X_0)_i .$$
By inequalities (3) one gets that
$\d(A_1, \dots ,A_{n-3}, X_0,X_0)_i =0$ for all $i$.
Hence $X_0$ vanishes by the induction hypothesis. This proves that the form $B$
is non-degenerate.

It remains to compute the signature of $B$. Remind that $B$
depends on positive definite matrices $A_1, \dots , A_{n-2}$. The
space of positive definite matrices is connected (indeed if $A$
and $B$ are positive definite then $tA+(1-t)B$ is positive
definite for $0\leq t \leq 1$). The signature of a family of
non-degenerate quadratic forms cannot jump. Hence it is constant.
Thus we can assume that $A_1= \dots =A_{n-2}=Id$. As in the case
$n=2$ it is sufficient to check that if  $X \ne 0$ satisfies
$B(X,Id)=0$ then $B(X,X) < 0$. Again we can assume that $X$ is
diagonal,
$$X= \left[ \begin{array}{ccc}
x_1& &0 \\
 &\ddots& \\
0& &x_n
\end{array}
\right] .$$
The condition $B(X,Id)=0$ means that $\sum_{i=1}^n x_i=0$.
Also it is easy to see that $$\kappa \cdot B(X,X)= 2 \sum_{i<j} x_i x_j ,$$
where $\kappa$ is a positive normalization constant.
But $$2 \sum_{i<j} x_i x_j =(\sum_i x_i)^2 -\sum_i x_i^2=
-\sum x_i ^2 <0.$$ The theorem is proved.
\qed

We will need also the following result.
\begin{theorem}
(i) The function $X\mapsto \log (detX)$ is concave on the cone of
positive definite hyperhermitian matrices, namely if $A,\,B\geq 0$
and $0\leq t\leq 1$ then
$$\log(det(tA+ (1-t)B))\geq t\log(detA) +(1-t)\log(detB).$$

(ii) The function $X\mapsto (\det X)^{\frac{1}{n}}$ is concave on the cone of the positive definite
hyperhermitian matrices.

(iii) If $A,\,B\geq 0$ then $$det(A+B)\geq detA +detB.$$
\end{theorem}
{\bf Proof.}
Note that we may assume that $A=I$ and $B$ is real diagonal. Both
results follow from the (known) real case. \qed


\subsection{ Dieudonn\'e determinant.}
We will remind the construction of the Dieudonn\'e determinant
referring for the details and proofs to \cite{artin}. Also we will
prove some properties of it which will be used in the subsequent
sections of the paper. Intuitively the Dieudonn\'e determinant of
an arbitrary quaternionic matrix has almost the same algebraic and
analytic properties as the {\itshape absolute value} of the usual
determinant of real or complex matrices. First let us discuss
purely algebraic construction.

Let $F$ be an infinite field, not necessarily commutative. Let
$M_n(F)$ denote the ring of $n \times n$-matrices with
coefficients in $F$. Let $GL_n(F)$ denote the group of invertible
$n \times n$-matrices. By an {\itshape elementary} matrix one
calls a matrix which has units on the diagonal and at most one
non-zero element out of the diagonal. Let $E_n$ denote the
subgroup of $GL_n(F)$ generated by all elementary matrices. Set
also $F^*_{ab} := F^*/[F^*,F^*]$ the abelinization of the
multiplicative group of $F$ (here $F^*$ denotes the multiplicative
group of $F$, and $[F^*,F^*]$ denotes its commutator subgroup).

\begin{theorem}[Dieudonn\'e]
Let $n\geq 2$. The group $E_n$ is normal subgroup of $GL_n(F)$.
For the quotient-group $GL_n(F) /E_n$ there exists a natural
isomorphism $D:GL_n(F) /E_n \str F^*_{ab}$.

This isomorphism $D$ is uniquely defined by the property that for
any invertible diagonal matrix $X= \left[ \begin{array}{ccc}
x_1& &0 \\
 &\ddots& \\
0& &x_n
\end{array}
\right] $,
$D(X)= \prod _i x_i \mbox{ mod } [F^*,F^*]$.
\end{theorem}

\begin{definition}[Dieudonn\'e determinant]
The Dieudonn\'e determinant
is a map $$D:M_n(F) \str F^*_{ab} \cup \{0\}$$
defined as follows: if $X$ is an invertible matrix
then $D(X)$ is as in Theorem 1.2.1; if $X$ is not
invertible then $D(X):=0$.
\end{definition}

Note also that it is convenient to define the Dieudonn\'e
determinant of elements of $F$, i.e. $1 \times 1$-matrices, as
$D(0)=0$ and for $x\ne 0$ as $D(x):=x \mbox{ mod } [F^*,F^*]$.

Let us state some basic general properties of the
Dieudonn\'e determinant. For the proofs we again refer to \cite{artin}.

\begin{theorem}

(i) $D(Id)=1$.

(ii) For $X,Y \in M_n(F)$
$$D(XY)= D(X) D(Y).$$

(iii) For any block-matrix $A= \left[ \begin{array}{cc}
X&0\\
0&Y
\end{array}
\right] $
with $X, \, Y$ being square matrices
$$D(A)= D(X) D(Y).$$

(iv) If one interchanges two rows or two columns
of the matrix then the Dieudonn\'e determinant
is multiplied by $-1 \mbox{ mod } [F^*, F^*]$.
\end{theorem}

Now let us consider in more details the case of quaternionic field
$F= \HH$. The commutator subgroup $[\HH ^*,\HH^*]$ coincides with
the subgroup of quaternions of absolute value 1. Thus we can
identify $\hab$ with the multiplicative group $\rp$ by
$$q \mbox{ mod } [\HH ^*,\HH^*] \leftrightarrow |q| :=\sqrt{q \bar q}.$$
So in the quaternionic case the Dieudonn\'e determinant maps
$$D:M_n(\HH) \str \rnn.$$
In the rest of the paper we will denote by $D(X)$ the Dieudonn\'e
determinant of a quaternionic matrix $X$, and by $det(X)$ the
Moore determinant of a {\itshape hyperhermitian} matrix $X$.


\begin{theorem}
(i) For any complex $n \times n$-matrix $X$ considered as
quaternionic matrix the Dieudonn\'e determinant $D(X)$ is equal to
the absolute value of the usual determinant of $X$.

(ii) Let $X$ be a quaternionic hyperhermitian $n \times n$-matrix.
Then its Dieudonn\'e determinant $D(X)$ is equal to the absolute
value of its Moore determinant $|det(X)|$.

(iii) For any $X$
$$D(X^*)=D(X),$$
where $X^*$ denotes
the transposed and quaternionic conjugate 
matrices respectively.\footnote{Added in Sept 2024: in the contrary to the commutative case, the equality $D(X^t)=D(X)$ is not satisfied in general, where $X^t$ is the transposed matrix of $X$; see  Section 4.2 in
	C.-Y. Lin and C.-F. Yu, Dieudonn\'e’s determinants and structure of
	general linear groups over division rings revisited, Bull. Inst. Math.
	Acad. Sin. (N.S.) 16 (2021), no. 1, 21–47.}
\end{theorem}

For any $n \times n$-matrix $X$ and any subsets $I, J \subset \{1,
\dots ,n\}$ let us denote by $M_{I,J}(X)$ the matrix obtained from
$X$ by deleting the rows with indexes  in $I$ and columns with
indexes in $J$. The following result is a weakened version of
usual formula of the decomposition of the determinant with respect
to a row. Note that this result is satisfied for the absolute
value of complex matrices.


\begin{theorem}
Let $A= \left[ \begin{array}{ccc}
a_{11}&\dots&a_{1n}\\
\multicolumn{3}{c}{\dotfill} \\
a_{n1}&\dots&a_{nn}
\end{array}
\right] $ be a quaternionic matrix.
Then $$D(A) \leq \sum_{i=1}^n |a_{1i}| D(M_{1i}).$$
Similar inequalities hold for any other row or column.
\end{theorem}

{\bf Proof.} From Theorem 1.2.3 it follows that
$$D \left( \left[ \begin{array}{c|ccc}
a&0&\dots&0\\ \hline
\ast& && \\
\vdots&&B& \\
\ast& &&
\end{array}
\right] \right) =|a|D(B).$$ Hence to prove the statement it is
sufficient to show that the Dieudonn\'e determinant is subadditive
with respect to the first row; namely if the matrices $A,A', A''$
are such that the first row of $A$ is the sum of first rows of
$A'$ and $A''$ and all the other rows are the same, then $D(A)
\leq D(A') +D(A'')$. But the Dieudonn\'e determinant has the
following property over arbitrary (non-commutative) field $F$
(\cite{artin}, Thm. 4.5):
$$D(A) \subset D(A') + D(A'') ,$$
where the inclusion and addition are understood in the sense of
conjugacy classes modulo $[F^*,F^*]$. But under our
identification of $\HH^*_{ab}$ with $\RR _{>0}$ the last inclusion
implies the desired inequality $D(A) \leq D(A') +D(A'')$. \qed

The next two propositions will be used in the sequel. It will be convenient
to introduce the following notation. Set $M'_{IJ}(A):=
M_{ \{1, \dots ,n\}-I,\{1, \dots ,n\}-J}(A)$, i.e.
it denotes the minor which stays on the intersection
of the rows with indexes from $I$ and columns with indexes from $J$.
\begin{proposition}
Let $A$ be hyperhermitian non-negative definite $n \times n$- matrix.
Fix an integer $k$, $1 \leq k \leq n$ and two subsets
$I,J \subset \{1, \dots ,n\}$ of cardinality $k$. Then
$$2D(M'_{IJ}(A)) \leq D(M'_{II}(A)) + D(M'_{JJ}(A)) .$$
\end{proposition}

{\bf Proof.} For simplicity of the notation
and without loss of generality we may assume that
$I\cup J=\{1, \dots ,n \}$, $I=\{1, \dots ,k\}$, and $J=\{n-k+1, \dots ,n\}$.

First let us reduce to the case $I\cap J=\emptyset$. We have
$$A= \left[ \begin{array}{c|c|c}
\ast&\ast&\ast\\ \hline
\ast&\mpij &\ast\\ \hline
\ast&\ast&\ast
\end{array}
\right].$$ For generic matrix $A$ the (hyperhermitian) minor
$\mpij$ is invertible. Then by  Claim 1.1.7 one can choose an
invertible matrix $U_0$ such that $U_0^* \mpij U_0 =Id$. Let $U=
\left[ \begin{array}{c|c|c} Id&0&0\\ \hline 0&U_0&0\\ \hline
0&0&Id
\end{array}
\right] .$  Consider matrix $A_1:=U^*AU$.
Clearly $D(\miin) =D(\mii) D(U)^2$, and similarly
for $M'_{IJ}$ and $M'_{JJ}$. Hence replacing $A$ by $A_1$
we may assume that $\mpij =Id$. Thus $A$ has the form
$$A= \left[ \begin{array}{c|c|c}
\ast&X&\ast \\ \hline
X^*&Id&Y^* \\ \hline
\ast&Y&\ast
\end{array}
\right].$$
Set $V= \left[ \begin{array}{c|c|c}
Id&-X&0\\ \hline
0&Id&0 \\ \hline
0&-Y&Id
\end{array}
\right] $.
Consider
$$A_2:= VAV^* =\left[ \begin{array}{c|c|c}
P&0&R \\ \hline
0&Id&0 \\ \hline
R^*&0&Q
\end{array}
\right] .$$
Here $P$ and $Q$ are hyperhermitian matrices.
Then $A_2$ has the same Dieudonn\'e
determinants of the minors $M'_{II}, \, M'_{IJ}, \, M'_{JJ}$
as $A$. Hence we may replace $A$ by $A_2$, and we will denote it
by the same letter $A$.
Then $\mii =\left[ \begin{array}{cc}
P&0\\
0&Id
\end{array}
\right]$,
$\mij =\left[ \begin{array}{cc}
0&R\\
Id&0
\end{array}
\right]$,
$\mjj =\left[ \begin{array}{cc}
Id&0\\
0&Q
\end{array}
\right]$.
So one has to show that
$$2 D(R) \leq D(P) + D(Q).$$
This inequality  is the statement of the proposition for the matrix
$\tilde A:= \left[ \begin{array}{cc}
P&R\\
R^*&Q
\end{array}
\right]$ which is also hyperhermitian
and positive definite since $A$ is.
Replacing $\tilde A$  by the matrix
$\left[ \begin{array}{cc}
U_1&0\\
0&U_2
\end{array}
\right] \tilde A
\left[ \begin{array}{cc}
U_1&0\\
0&U_2
\end{array}
\right] ^* $
with $U_1, \, U_2 \in Sp(k)$ one can assume that
the matrices $P$ and $Q$ are diagonal.

Fix now some $U,V \in Sp(k)$ (the choice of them will be clear later).
Let $T:= \left[ \begin{array}{cc}
P^{1/2}UP^{-1/2}&0\\
0&Q^{1/2}VQ^{-1/2}
\end{array}
\right]$. Then
$$T \tilde A T^* = \left[ \begin{array}{cc}
P&R_1 \\
R_1^*&Q
\end{array}
\right],$$
where $R_1=P^{1/2}U(P^{-1/2}R Q^{-1/2})V^* Q^{1/2}$.
Note that $D(R_1)=D(R)$. Since $P$ and $Q$
are diagonal, by a choice of $U,V \in Sp(k)$
one can make the matrix $R_1$ diagonal.

Finally we are reduced to the hyperhermitian
non-negative definite matrix $A$ of the form
$A= \left[ \begin{array}{ccc|ccc}
\lam_1& &0&\nu_1& &0\\
 &\ddots& &&\ddots & \\
0& &\lam_k&0& &\nu_k\\ \hline
\bar \nu_1& &0&\mu_1& &0\\
 &\ddots& &&\ddots & \\
0& &\bar \nu_k &0& &\mu_k
\end{array}
\right]$. We have to show that
$$2 \prod_1^k |\nu _i| \leq \prod_1^k |\lam _i| +\prod_1^k |\mu _i|.$$
Consider the $2\times 2$- matrix
 $\left[ \begin{array}{cc}
\lam_i &\nu _i \\
\bar \nu _i&\mu _i
\end{array}
\right]$
which is clearly non-negative definite. Take a vector
$\left( \begin{array}{c}
1\\
t\cdot q
\end{array}
\right)$
for any $t \in \RR $ and any quaternion $q$ of norm 1.
Applying that matrix to this vector we get
$$ \lam_i +t^2 \mu_i +2t Re(\nu_i q) \geq 0 .$$
Hence $|\nu_i| \leq \sqrt{\lam _i \mu _i}$.
Then $$2 \prod _i |\nu_i| \leq 2 \sqrt{\prod _i |\lam_i| \cdot
\prod _i |\mu_i| }\leq \prod _i |\lam _i| +\prod _i |\mu_i|.$$
\qed
\begin{proposition}
Let $A=(a_{ij}), B$ be $n\times n$-hyperhermitian matrices. Then
the mixed discriminant satisfies
$$ |det(A, \underset{n-1 \mbox{ times}} {\underbrace{B, \dots,B}})|
\leq c_n \cdot max _{i,j}|a_{ij}| \cdot (\sum_{|I|,|J| =n-1} D(M'_{IJ}(B))),$$
where $c_n$ is a  constant depending on $n$ only.
\end{proposition}

{\bf Proof.} Since $det(A,B, \dots,B)$ is linear in $A$
it is sufficient to prove the inequality in the
following two cases:

1) $A= \left[ \begin{array}{cccc}
1& & & 0\\
 &0& & \\
 & &\ddots& \\
 0& & &0
\end{array}
\right]$;
2) $A= \left[ \begin{array}{cc|ccc}
0&q& && 0\\
\bar q&0& & & \\ \hline
& &0& &\\
 & & & \ddots & \\
0 & & & &0
\end{array}
\right]$.

 The first case follows from Proposition 1.1.11.
Let us consider the second case. Replacing
$A$ by the matrix
$$\left[ \begin{array}{cccc}
\frac{\bar q}{|q|}& & &0\\
&1& &\\
&&\ddots& \\
0& & &1
\end{array}
\right] A \left[ \begin{array}{cccc}
\frac{ q}{|q|}& & &0\\
&1& &\\
&&\ddots& \\
0& & &1
\end{array}
\right]=
\left[ \begin{array}{cc|ccc}
0&|q|&& &0\\
|q|&0&& &\\ \hline
&&&\ddots& \\
0& && &0
\end{array}
\right]$$
we can assume that
$A= \left[ \begin{array}{cc|ccc}
0&1& && 0\\
1&0&& & \\ \hline
&&0& &\\
 & && \ddots & \\
0 & & &&0
\end{array}
\right]$.
Let $B= \left[ \begin{array}{c|c}
P&R\\ \hline
R^*&Q
\end{array}
\right]$. Here $P$ and $Q$ are hyperhermitian
matrices of sizes $2\times 2$ and $(n-2)\times (n-2)$ respectively.

{\bf Claim.} $$det(A,B, \dots ,B) \leq M'_{\{2,3, \dots ,n\}, \{1,3, \dots,n\}}(B) +
M'_{\{1,3, \dots ,n\}, \{2,3, \dots,n\}}(B).$$
 It remains to prove this claim. We may also assume that $Q$ is invertible.
Set $S:= \left[ \begin{array}{c|c}
Id& -R Q^{-1}\\ \hline
0&Id
\end{array}
\right]$. Consider $$B_1:=SBS^*= \left[ \begin{array}{cc}
\ast&0\\
0&\ast
\end{array}
\right].$$
Note also that $S^*AS=A$.
It is easy to see that
$$M'_{\{2,3, \dots ,n\}, \{1,3, \dots,n\}}(B_1)=
M'_{\{2,3, \dots ,n\}, \{1,3, \dots,n\}}(B) \mbox{ and }$$
$$M'_{\{1,3, \dots ,n\}, \{2,3, \dots,n\}}(B_1)=
M'_{\{1,3, \dots ,n\}, \{2,3, \dots,n\}}(B).$$
Hence it is sufficient to prove the claim under
assumption $R=0$, i.e.
$B= \left[ \begin{array}{cc}
P&0\\
0&Q
\end{array}
\right]$.
Then clearly $$det(A,B, \dots ,B)=
det \left( \left[ \begin{array}{cc}
0&1\\
1&0
\end{array}
\right], P \right) \cdot det Q.$$
If $P=\left[ \begin{array}{cc}
b_{11}&b_{12}\\
b_{21}&b_{22}
\end{array}
\right]$ then
$$|det \left(\left[ \begin{array}{cc}
0&1\\
1&0
\end{array}
\right], P \right)|= Re (b_{12}) \leq |b_{11}|+|b_{22}|,$$
where the last inequality follows from Proposition 1.2.6.
Proposition 1.2.7 is proved.
\qed

From Propositions 1.2.6 and 1.2.7 one can easily deduce
\begin{proposition}
Let $A= (a_{ij})$ be a hyperhermitian matrix
and $B_1, \dots , B_{n-1}$ be non-negative
definite hyperhermitian matrices. Then
\begin{eqnarray*}
|det(A,B_1, \dots , B_{n-1})| \leq\\
c_n \cdot max_{i,j} |a_{ij}| \cdot \sum_{|I|=n-1}\sum_{1 \leq i_1, \dots ,i_{n-1} \leq n-1}
det(M_{II}'(B_{i_1}), \dots,M_{II}'( B_{i_{n-1}})) ,
\end{eqnarray*}
where $c_n$ is a constant depending on $n$ only.
\end{proposition}


\section{Plurisubharmonic functions of quaternionic
variables.}  \setcounter{theorem}{0} In this part we will develop
a basic theory of \psh functions of quaternionic variables.

\subsection{Main notions.} First let us remind few standard notions. Below
$\Omega$ will denote an open domain. As usual we will denote by
$C^k(\Omega)$ the class of $k$ times continuously differentiable
functions on $\Omega$, and by $C^k_0(\Omega)$ the class of $k$
times continuously differentiable functions on $\Omega$ with
compact support. We will also denote by $L^\infty(\Omega)$ (resp. $L^\infty _{loc}(\Omega)$ )
 the class of bounded (resp. locally bounded) measurable functions on $\Omega$.
\begin{definition}
A real valued  function
$f: \Omega \subset \RR ^m \str \RR$ is called {\itshape subharmonic}
if

(a) $f$ is upper semi-continuous, i.e. $f(x_0)\geq \underset{x\str x_0}{\limsup f(x)}$ for any $x_0\in \Omega$;

(b) $f(x_0)\leq \int_{S(x_0,r)} f(x) d\sigma$ for any point $x_0$ and for any sufficiently small $r>0$. Here
$S(x_0,r)$ denotes the sphere of radius $r$ with center at $x_0$, and $\sigma$ is the Lebesgue measure on it normalized
by one.

\end{definition}
\begin{definition}
A real valued continuous function
$$f: \Omega \subset \RR ^n \str \RR$$ is called {\itshape
convex} if its restriction to any (real) line is subharmonic.
\end{definition}
\begin{definition}
A real valued  function
$$f: \Omega \subset \CC ^n \str \RR$$ is called {\itshape
\psh} if it is upper semi-continuous and its restriction to any {\itshape complex} line is subharmonic.
\end{definition}
Now let us introduce a new definition.
\begin{definition}
A real valued function
$$f: \Omega \subset \HH ^n \str \RR$$
is called {\itshape quaternionic \psh} if it is upper
semi-continuous and its restriction to any right
 {\itshape quaternionic}
line is subharmonic.
\end{definition}

It is easy to see that any (quaternionic) \psh function is subharmonic.

\begin{example}
1) Any convex function on $\HH ^n$ is quaternionic \psh.

2) Fix on $\HH ^n$ one of the complex structures compatible with
the quaternionic structure; say, let us fix $i$. Let $f$ be a
 \psh function with respect to this complex
structure in the sense of Definition 2.1.3. It is easy to see that
$f$ is \psh in the quaternionic sense.
\end{example}

Let $q$ be a quaternionic coordinate,
$$q=t+ix+jy +kz ,$$
where $t,x,y,z$ are real numbers. Consider the following operators
defined on the class of smooth $\HH$-valued functions of the
variable $q\in \HH$:
$$\db f:=\frac{\partial f}{\partial  t}  +
i \frac{\partial f}{\partial  x} +
j \frac{\partial f}{\partial  y} +
k \frac{\partial f}{\partial  z}  , \mbox { and }$$
$$\dq f:=\overline{ \db \bar f}=
\frac{\partial f}{\partial  t}  -
 \frac{\partial f}{\partial x}  i-
 \frac{\partial f}{\partial  y} j-
\frac{\partial f}{\partial  z}  k.$$ Note that $\db$ is called
sometimes Cauchy-Riemann-Moisil-Fueter operator, and sometimes
Dirac-Weyl operator (see the introduction). It is easy to see that
$\db$ and $\dq$ commute, and if $f$ is a {\itshape real valued}
function then $$\db \dq f= \Delta f = (\frac{\partial ^2
}{\partial t ^2}+ \frac{\partial ^2 }{\partial x ^2}+
\frac{\partial ^2 }{\partial y ^2}+ \frac{\partial ^2 }{\partial z
^2})f.$$ For any real valued $C^2$- smooth function $f$ the matrix
$(\dfq)_{i,j=1}^{n}$ is obviously hyperhermitian. For brevity we
will use the following notation:
$$det(f_1, \dots, f_n):=
det \left( (\frac{\partial ^2 f_1}{\partial \bar q_i \partial  q_j}), \dots,
(\frac{\partial ^2 f_n}{\partial \bar q_i \partial  q_j}) \right) ,$$
where $det$ denotes the mixed discriminant of hyperhermitian matrices
(see Definition 1.1.14). Note also that the operators $\frac{\partial}{\partial q_i}$ and
$\frac{\partial}{\partial \bar q_j}$ commute.
One can easily check the following identities.

{\bf Claim.} {\itshape (i) Let $f:\HH ^n \str \HH$ be a smooth
function. Then for any $\HH$-linear transformation $A$ of $\HH ^n$
(as right $\HH $-vector space) one has the identities
$$ \left( \frac {\partial ^2 f(Aq)}{\partial \bar q_i \partial q_j} \right)
=A^* \left(\frac {\partial ^2 f}{\partial \bar q_i \partial q_j}(Aq) \right)A .$$

(ii) If, in addition, $f$ is real valued then for any $\HH$-linear
transformation $A$ of $\HH ^n$ and any quaternion $a$ with $|a|=1$
$$ \left( \frac {\partial ^2 f(A(q \cdot a))}{\partial \bar q_i \partial q_j} \right)
=A^* \left(\frac {\partial ^2 f}{\partial \bar q_i \partial q_j}(A(q\cdot a)) \right) A .$$
}

\begin{proposition}
A real valued twice continuously
differentiable function $f$ on the domain $\Omega \subset \HH ^n$
is quaternionic
plurisubharmonic if and only if at every point $q \in \Omega$
the matrix $(\frac{\partial ^2 f}{\partial \bar q_i \partial  q_j})(q)$
is non-negative definite.
\end{proposition}
The proof of this proposition is straightforward. The following
lemma will be useful in the sequel.

\begin{lemma}
Let $f_0, f_1, \dots, f_n$ be real valued compactly supported sufficiently
smooth functions on $\HH ^n$.
The $(n+1)$- linear functional
$$L(f_0, f_1, \dots ,f_n) :=
\int _{\HH ^n} f_0(q) \cdot det(f_1, \dots, f_n) (q) dq $$
is symmetric with respect to all $f_0, f_1, \dots, f_n$.
\end{lemma}

{\bf Proof.} Note that $L$ is symmetric with respect to the
last $n$ arguments. Thus it is sufficient to check that
\begin{equation}
L(f_0,f_1,f_2, \dots ,f_n)= L(f_1,f_0,f_2, \dots ,f_n)
\end{equation}
for any smooth compactly supported functions $f_0, f_1, \dots
,f_n$. Both sides of (4) make sense if $f_0$ is a generalized
function. Since linear combinations of delta-functions of points
$\delta _{q}$ are dense in the space of all the generalized
functions it is sufficient to prove (4) for $f_0= \delta _{0}$,
namely
\begin{equation}
(det(f_1, \dots , f_n))|_{q=0}=
\ih f_1 (q) det ( \delta _0, f_2, \dots ,  f_n).
\end{equation}
Clearly the right hand side in equation (5) depends
only on derivatives at 0 of $f_1, \dots ,f_n$   up to order 2.
Consider the terms of the Taylor series of $f_1$ at 0:
$$f_1(q) =g(q)+h(q)+ O(|q|^3),$$
where $g$ is a polynomial of degree one, and  $h$
is a quadratic term. So it is sufficient to prove the following two statements:

Case 1. \begin{equation}
L(h, \delta_0, f_2, \dots ,f_n)= det(h, f_2, \dots ,f_n)|_{q=0}
\end{equation}
for any smooth compactly supported function $h$ which is equal to
a homogeneous polynomial of degree 2
 in a neighborhood of 0, and for any smooth compactly
supported functions $f_2, \dots ,f_n$.

Case 2. \begin{equation}
L(g, \delta _0,f_2, \dots ,f_n)=0
\end{equation}
for any smooth compactly supported function $g$ which is equal to
a polynomial of degree 1 in a neighborhood of 0, and for any
smooth compactly supported functions $f_2, \dots ,f_n$.

Let us consider Case 1.
If we write down the formula for $L(h,\delta _0, f_2, \dots ,f_n)$
as a polynomial in $\frac {\partial ^2 f_k}{\partial t_i
\partial t_j}, \, \frac {\partial ^2 f_k}{\partial t_i
\partial x_j} $ etc. and in
 $\frac {\partial ^2 \delta _0}{\partial t_i
\partial t_j}, \, \frac {\partial ^2 \delta _0}{\partial t_i
\partial x_j}$ etc. then we see that
the derivatives of
$\delta _0$ enter at each monomial only once
 because of linearity of $L$ with respect to each argument. For example
consider a monomial containing $\frac {\partial ^2 \delta _0}{\partial t_i
\partial t_j}$. Let it be $ \ih h \cdot
\frac {\partial ^2 \delta _0}{\partial t_i \partial t_j}
\cdot \partial ^2 f_2 \cdot \dots \cdot \partial ^2 f_n$,
where $\partial ^2 f_k$ denotes certain partial derivative
of order 2 of $f_k$.
But $$\ih h \cdot
\frac {\partial ^2 \delta _0}{\partial t_i \partial t_j}
\cdot \partial ^2 f_2 \cdot \dots \cdot \partial ^2 f_n =
\frac {\partial ^2}{\partial t_i \partial t_j}
(h \cdot \partial ^2 f_2 \cdot \dots \cdot \partial ^2 f_n)|_{q=0}=$$
$$ \frac {\partial ^2 h}{\partial t_i \partial t_j}(0)
\cdot  \partial ^2 f_2(0) \cdot \dots \cdot \partial ^2 f_n (0) ,$$
where the last equality is satisfied since
the first derivatives of $h$ at 0 vanish.
Thus in each monomial the term $h \cdot
 \frac {\partial ^2 \delta _0}{\partial t_i \partial t_j}$
is just replaced by $\frac {\partial ^2 h}{\partial t_i \partial t_j}(0)$.
Hence the final expression is $det(h,f_2, \dots ,f_n) |_{q=0}$.
This proves the first case.

Let us prove Case 2. It is convenient to prove a more general statement.

{\bf Claim.} {\itshape Let $U$ be a fixed neighborhood of the
origin 0. Let $g$ be any smooth compactly supported function
which is equal to a polynomial of degree 1 inside $U$. Let $f_1$
be a generalized function with support contained in $U$. Let $f_2,
\dots ,f_n$ be smooth compactly supported functions.

Then $$ \ih g \det(f_1, f_2, \dots , f_n) =0.$$
}
The  proof of the claim will be by induction in $n$.
If $n=1$ then using selfadjointness of the Laplacian one gets:
$$ \int _{\HH} g \Delta f_1= \int _{\HH} \Delta g \cdot f_1=
\int _{U} \Delta g \cdot f_1 =0.$$ Assume that $n>1$. It is well
known (see Appendix) that the linear combinations of
delta-functions of quaternionic hyperplanes are dense in the space
of all generalized functions (this fact is equivalent to the
injectivity of the Radon transform with respect to quaternionic
hyperplanes). Hence it is sufficient to prove the claim for $f_1
=\delta _{L}$, where $L$ is the hyperplane $\{q_1=0\}$.

Since $\delta _{L}$ is invariant with respect
to translations in directions $q_2, \dots , q_n$ then
$\frac {\partial ^2 \delta _{L}}{\partial \bar q_i \partial  q_j} =0$ unless $i=j=1$.
Thus
$$\left(\frac {\partial ^2 \delta _{L}}{\partial \bar q_i \partial q_j}\right) =
\left[ \begin{array}{c|ccc}
\Delta _1 \delta _{L}&0& \dots &0 \\ \hline
0&&&\\
\vdots&&0&\\
0&&&
\end{array} \right],$$
where $\Delta _1$ denotes the Laplacian with respect to
the first coordinate: $\Delta _1= \frac{\partial ^2}{\partial t_1^2}+
\frac{\partial ^2}{\partial x_1^2}+
\frac{\partial ^2}{\partial y_1^2}+
\frac{\partial ^2}{\partial z_1^2}$.
Using Proposition 1.1.11 it is easy to see that
$$c \cdot det (\delta _{L}, f_2, \dots ,f_n)=
\Delta _1 \delta _{L} \cdot det(B_2, \dots ,B_n),$$
where $c$ is a positive
normalizing constant, and $B_k$ denotes the $(n-1) \times (n-1)$- matrix
$(det \frac{\partial ^2 f_k}{\partial \bar q_i
\partial  q_j})_{i,j=2}^{n}$.
Then $$\ih g \cdot det (\delta _{L},f_2, \dots ,f_n)= \ih g \cdot
\Delta _1 \delta _{L} \cdot det(B_2, \dots ,B_n).$$ Clearly the
last expression depends only on the 2-jets of $g, f_2, \dots ,f_n$
in the direction $q_1$. Thus we may assume that the functions
$f_k$ are of the form
$$f_k(q_1,q_2, \dots ,q_n)=p_k (q_1)
\cdot f'_k (q_2, \dots ,q_n),$$
where $p_k(q_1)$ are polynomials
(of degree at most 2) depending only on $t_1,x_1,y_1,z_1$,
and $f_k '$ are smooth compactly supported functions
 depending only on $q_2, \dots ,q_n$.

Since $deg \, g \leq 1$ we may assume (by linearity)
that either $g(q_1,q_2, \dots , q_n)= g(q_1)$
or $g(q_1,q_2, \dots , q_n)= g(q_2, \dots ,q_n)$.
In the first case
$$\ih g \cdot \Delta _1 \delta _{L} \cdot det(B_2, \dots ,B_n)=$$
$$ \Delta _1 \left( g(q_1) \cdot p_2(q_1) \cdot \dots \cdot p_n(q_1)\right) |_{q_1=0}
\cdot \int_{L} det(B'_2, \dots ,B'_n),$$ where $B'_k$ denotes the
matrix $(\frac{\partial ^2 f'_k}{\partial \bar q_i \partial
q_j})_{i,j=2}^{n}$. The last integral vanishes by the induction
assumption.

Now consider the second case $g(q_1,q_2, \dots , q_n)= g(q_2, \dots ,q_n)$.
We have
$$\ih g \cdot \Delta _1 \delta _{L} \cdot det(B_2, \dots ,B_n)=
\Delta _1 ( p_2 \dots p_n)|_{q_1=0}
\int_{L}g\cdot det(B_2, \dots ,B_n).$$
Again the last expression vanishes by the induction assumption.
Thus our claim, and hence Proposition 2.1.6, are proved.
\qed

The next result is again a quaternionic analogue of
the corresponding property of convex functions and complex
plurisubharmonic functions. We adopt the arguments of
Chern-Levine-Nirenberg \cite{chern-levine-nirenberg} and
 Rauch-Taylor
\cite{rauch-taylor}
(see also \cite{aubin}).


\begin{proposition}
Let $\Omega \subset \HH ^n$ be an open domain. Assume that a
sequence $\{f_N \}$ of twice continuously differentiable
quaternionic plurisubharmonic functions converges uniformly on compact subsets
 to a twice continuously differentiable
function $f$. Then $f$ is also quaternionic plurisubharmonic, and
for every continuous function $\phi$ with compact support in
$\Omega$
$$\int_{\Omega} \phi \cdot det(\frac{\partial ^2 f_N}
{\partial \bar q_i \partial  q_j})\str
\int_{\Omega} \phi \cdot det(\frac{\partial ^2 f}
{\partial \bar q_i \partial  q_j}) \mbox{ as } N\str \infty.$$
\end{proposition}

We will need a lemma. But first let us introduce a notation. For
subsets $I,J \subset \{1, \dots ,n\}$ and a function $g$ let us
denote by $M'_{IJ}(g)$ the matrix which stays on the intersection
of rows with indexes from $I$ and columns with indexes from $J$ in
the matrix $det(\frac{\partial ^2 g}{\partial\bar q_i \partial
q_j})_{i,j=1}^{n}$. Also for a set $U$ and a function $g$ defined
on it let us denote by $||g|| _{L^\infty (U)} :=sup _{q\in U}
|g(q)|$, and by $||g||_{C^k(U)}$ the maximum of
$L^\infty(U)$-norms of all partial derivatives of $g$ up to order
$k$.
 Below we will
denote for brevity
$det(\frac{\partial ^2 g}{\partial\bar q_i \partial q_j})_{i,j=1}^{n}$
by $det(g)$.

\begin{lemma}
Let $I, \, J$ be subsets of $\{1, \dots ,n\}$ of cardinality $k$.
Let $f\in L^\infty _{loc}(\Omega)$, and let $ g$ be a
twice continuously differentiable quaternionic \psh function on a
domain $\Omega \subset \HH ^n$. Let $K$ be a compact subset of
$\Omega$, and let $U$ be a compact neighborhood of $K$ in $\Omega$. Then
$$ |\int _K f \cdot D(M'_{IJ}(g))| \leq C(U) ||f||_{L^\infty(K)} ||g||_{L^\infty(U)} ^k ,$$
where $C(U)$ is a constant depending on $U$ only.
\end{lemma}

{\bf Proof.} Since $g$ is \psh , Proposition 1.2.6 implies the
estimate $D(M'_{IJ}(g)) \leq D(M'_{II}(g)) +D(M'_{JJ}(g))$. Hence
$$|\int _K f \cdot D(M'_{IJ}(g))| \leq ||f||_{L^\infty(K)} \cdot \int
_K (D(M'_{II}(g)) +D(M'_{JJ}(g))) .$$ It remains to prove that for
any subset $I\subset \{1, \dots, n\}$ of cardinality $k$
$$\int _K det(M'_{II}(g)) \leq C(U) ||g||_{L^\infty(U)} ^k .$$
Let us prove this inequality by induction in $k$. For $k=0$
the statement is trivial. Assume that $k>0$.
Let us fix a non-negative function $\gamma \in C^{\infty}_0(\Omega)$
such that $\gamma |_K \equiv 1$ and $\gamma$ vanishes
on $\Omega -U$. Then using Lemma 2.1.7
$$\int _K det M'_{II}(g) \leq \ih \gamma \cdot det M'_{II}(g)
= \ih g \cdot det_I (\gamma, g, \dots , g) , $$ where $det_ I$
denotes the mixed discriminant of matrices of order $|I|$. By
Proposition 1.2.7 the last expression is at most
$$ ||g||_{L^\infty(U)} \cdot ||\gamma||_{C^2(U)}
\sum _{S,T} \int _U D_{ST}(g), $$
where the sum extends over all subsets $S,\, T$ of $I$
of cardinality $k-1$. Again by Proposition 1.2.6
$$\int _U D_{ST}(g) \leq \int _U (det_S(g) +det_T (g)).$$
Now the estimate follows by the assumption of induction.
\qed

Now let us prove Proposition 2.1.8.
 First let us show that the
limit $f$ is plurisubharmonic. This is obvious since the
restriction of $f$ to any quaternionic line is subharmonic as the
uniform limit of subharmonic functions.

Let us prove the second part of Proposition 2.1.8. Let $K:=supp
\phi$. Fix $\eps
>0$, and a compact neighborhood $U$ of $K$. Let us choose a
function $\psi \in C^{\infty}_0(\Omega)$ such that $||\phi
-\psi||_{L^\infty(U)} \leq \eps$. We have $$|\int _K \psi (det
(f_N) -det (f)) - \int _K \phi (det (f_N) -det (f))| \leq $$
$$ |\int _K (\psi -\phi) det( f_N) | +
|\int _K (\psi -\phi) det(f)| \leq$$
$$ C(U) (||f_N||_{L^\infty(U)}^n +||f||_{L^\infty(U)}^n) \cdot \eps ,$$
where the last inequality follows from Lemma 2.1.9.
The last expression can be estimated for large $N$ by
$3C(U) ||f||_{L^\infty(U)}^n \cdot \eps$. Thus it is sufficient
to prove that
$$\int_{\Omega} \psi \cdot det(f_N)
\mbox{ tends to }
\int_{\Omega} \psi \cdot det(f) \mbox{ as } N\str \infty.$$
We have
$$ |\int _{\Omega} \psi \cdot (det(f_N) -det(f))| =|\int _{\Omega}
\sum_{i=0}^{n-1} \psi \cdot det(\underset{i \mbox{ times }}
{ \underbrace {f_N, \dots, f_N}}, f_N -f,
\underset {n-i-1 \mbox { times }} {\underbrace{
f, \dots ,f}})| = $$
$$|\sum_{i=0}^{n-1} \int _{\Omega} (f_N-f) det(\underset{i}
{\underbrace{f_N, \dots ,f_N}}, \psi,
\underset {n-i-1 } {\underbrace{
f, \dots ,f}})| \leq$$
\begin{equation}
 C \cdot ||f_N -f|| |_{L^\infty(\Omega)}
||\psi||_{C^2(\Omega)}
\sum_{i=0}^{n-1} \sum_{|I|=n-1}
\int _{supp \, \psi} det_I(\underset{i}
{\underbrace{f_N, \dots ,f_N}}, \underset {n-i-1 } {\underbrace{
f, \dots ,f}})
\end{equation}
by Proposition 1.2.8 (here we have used the fact that the
functions $f_N$ and $f$ are plurisubharmonic). Now let us estimate
the last expression.
$$\int _{supp \, \psi}det_I (f_N,\dots,f_N,f,\dots,f)
\leq \int _{supp \, \psi}det_I(\underset{n-1}
{\underbrace{f+f_N, \dots,f+f_N}})
\leq $$
$$C' ||f+f_N||_{L^\infty}^{n-1},$$
where the last inequality holds by Lemma 2.1.9.
Hence the expression (8) tends to 0 as $N\str \infty$.
This proves Proposition 2.1.8.
\qed

Now let us study {\itshape continuous} quaternionic \psh functions
which are not necessarily smooth. For every continuous \psh
function $f$ we will define a non-negative measure such that if
$f$ is smooth it coincides with $det(\dfq)$. To do it let us
observe first of all that any continuous \psh function $f$ on a
domain $\Omega \subset \HH^n$ can be approximated by $C^{\infty}$-
smooth \psh functions uniformly on  compact subsets of $\Omega$.
(To see it consider the convolution of $f$ with the delta-sequence
of non-negative $C^{\infty}$-smooth functions. Each such
convolution is infinitely smooth and \psh). The next theorem is
first main result of this section; it provides the definition of
the measure $det(\frac{\partial ^2 f}{\partial\bar q_i \partial
q_j})$ for any continuous \psh function $f$.

\begin{theorem}
Let $f$ be a continuous quaternionic \psh
function on a domain $\Omega$. Let $\{f_N\}$
be a sequence of twice continuously differentiable
\psh functions converging to $f$ uniformly
on compact subsets of $\Omega$. Then
$det(\frac{\partial ^2 f_N}{\partial\bar q_i \partial  q_j})$
weakly converges to a non-negative measure on $\Omega$.
This measure depends only on $f$ and not on the choice
of an approximating sequence $\{f_N\}$.
\end{theorem}

This measure will be denoted by $det(\dfq)$.

{\bf Proof.} By Lemma 2.1.9 one sees that for any compact subset
$K\subset \Omega$ the sequence of measures $det(\frac{\partial ^2
f_N}{\partial\bar q_i \partial  q_j})|_K$ is bounded. Thus it is
sufficient to show that for any continuous compactly supported
function $\phi$ the sequence $\int _{\Omega} \phi \cdot
det(\frac{\partial ^2 f_N}{\partial\bar q_i \partial  q_j})$ is a
Cauchy sequence. Let us fix $\eps >0$, and a function $\psi \in
C^\infty _0(\Omega)$ such that $||\phi -\psi||_{C^0(\Omega)} <
\eps$. Let us also fix an arbitrary compact subset $K\subset
\Omega$ and a compact neighborhood $U$ of $K$ in $\Omega$. As in
the proof of Proposition 2.1.8 we have
$$|\int _K \psi (det (f_N) -det (f_M)) -
\int _K \phi (det (f_N) -det (f_M))| \leq $$
$$ |\int _K (\psi -\phi) det( f_N) | +
|\int _K (\psi -\phi) det(f_M)| \leq$$
$$ C(U) (||f_N||_{C^0(U)}^n +||f_M||_{C^0(U)}^n) \cdot \eps ,$$
where the last inequality follows from Lemma 2.1.9. For large $M$
and $N$ the last expression can be estimated from above by
$3C(U) ||f||_{C^0(U)}^n \cdot \eps$.
Hence it is sufficient to prove that for any
function $\psi \in C^\infty _0 (\Omega)$
the sequence  $\int _{\Omega} \psi \cdot
det((\frac{\partial ^2 f_N}{\partial\bar q_i \partial  q_j}))$
is a Cauchy sequence. We have the following
estimate exactly as in the inequality (8) (with $f_M$
instead of $f$):
$$ |\int _{\Omega} \psi \cdot (det(f_N) -det(f_M))| \leq$$
$$ C \cdot ||f_N -f_M|| |_{C^0(supp \psi)}
||\psi||_{C^2(\Omega)}
\sum_{i=0}^{n-1}\sum_{|I|=n-1}
\int _{supp \, \psi} det_I(\underset{i}
{\underbrace{f_N, \dots ,f_N}}, \underset {n-i-1 } {\underbrace{
f_M, \dots ,f_M}}).$$
Again as in the proof of Proposition 2.1.8
we get
$$\int _{supp \, \psi} det_I(\underset{i}
{\underbrace{f_N, \dots ,f_N}}, \underset {n-i-1 } {\underbrace{
f_M, \dots ,f_M}}) \leq
C||f_N+f_M||^{n-1}_{C^0(supp \psi)} < C'.$$
This proves Theorem 2.1.10.
\qed

The second main result of this section is as
follows.
\begin{theorem}
Let $\{f_N\}$ be a sequence of continuous quaternionic \psh
functions in a domain $\Omega \subset \HH^n$. Assume that this
sequence converges uniformly on compact subsets to a function $f$.
Then $f$ is continuous quaternionic \psh function. Moreover the
sequence of measures $det(\frac{\partial ^2 f_N}{\partial\bar q_i
\partial  q_j})$ weakly converges to the measure $det(\dfq)$.
\end{theorem}

{\bf Proof.} The limit $f$ is a \psh function. Indeed
the restriction of $f$ to any quaternionic line is subharmonic
as a uniform limit of subharmonic functions.

Let us prove the second part of the statement.
 The functions $f_N$ can be approximated
uniformly on compact subsets as good as we wish
by smooth \psh functions $g_N$ such that
the sequence $g_N$ will converge uniformly on compact subsets to
$f$. Then the result follows from previous Theorem 2.1.10.
\qed
\subsection{The minimum principle.} In this subsection we prove the
following minimum principle.
\def\bfm{{\cal B} (\phi, f)}
\def\bo{\partial \Omega}
\def\co{\bar \Omega}
\def\vz{v_{\zeta}}
\def\fss{F^{**}}
\def\ggo{\Gamma_0}
\def\pt{\partial ^2}
\def\oe{\omega(\eps)}
\def\xom{\chi_{\Omega}}
\begin{theorem}
Let $\Omega$ be a bounded  open set in $\HH^n$. Let $u,\, v$ be
continuous functions on $\bar \Omega$ which are \psh in $\Omega$.
Assume that
$$det(\frac{\partial ^2 u}{\partial\bar q_i \partial  q_j})\leq
det(\frac{\partial ^2 v}{\partial\bar q_i \partial  q_j}) \mbox{
in } \Omega.$$ Then
$$min\{u(z)-v(z)|z\in \bar \Omega\} =min\{u(z)-v(z)|z\in \partial \Omega\}.$$
\end{theorem}
The exposition follows very closely to Section 3 of
\cite{bedford-taylor}. From now on we will denote for brevity the
matrix $\frac{\partial ^2 u}{\partial\bar q_i \partial q_j}$ by
$\pt u$.

\begin{proposition}
Let $\Omega$ be a bounded domain in $\HH ^n$ with smooth boundary,
and let $u,v\in C^2(\bar \Omega)$ be psh functions on $\Omega$. If
$u=v$ on $\partial \Omega $ and $u\geq v$ in $\Omega$, then
$$ \int _{\Omega} det (\partial^2 u) \leq \int _{\Omega} det (\partial^2 v).$$
\end{proposition}
{\bf Proof.} First we can write $\Omega=\{\rho<0\}$ with $\rho$
being a smooth function, $\bo =\{\rho =0\}$, and $\nabla
\rho|_{\bo}\ne 0$. We have
$$\int _{\Omega} (det (\pt u)-det(\pt v))=\sum_{i=0}^{n-1}\int _{\Omega}
det(\underset{i\mbox{ times}}{\underbrace{ u, \dots,  u}},
\underset{n-i-1\mbox{ times}}{\underbrace{ v, \dots,  v}},
 (u-v)).$$ Let us prove that each summand is non-positive. We
will need the following lemma.
\begin{lemma}
Let $\beta\in C^2(\co),\, \beta |_{\bo}\equiv 0$. Let $u_1,\dots,
u_{n-1}\in C^3(\co)$. Then
$$ \int_{\Omega} det( u_1,\dots,  u_{n-1}, \beta)=
-\int_{s\in\bo} det ( u_1|_{\tilde T_s},\dots,
u_{n-1}|_{\tilde T_s}) \frac{\partial \beta}{\partial \nu
(s)}ds,$$ where $T_s$ denotes the tangent space to $\bo$ at $s$,
$\tilde T_s$ denotes the quaternionic subspace $T_s\cap i\cdot
T_s\cap j\cdot T_s \cap k\cdot T_s$, $\nu (s)$ is the inner normal
to $\bo$, and $ds$ is the surface area measure.
\end{lemma}
Let us continue proving Proposition 2.2.2 assuming this lemma.
Since $u\geq v$ we can represent $u-v=\alpha\cdot \rho$, where
$\alpha\leq 0$. Using Lemma 2.2.3 we have
$$\int _{\Omega}
det(\underset{i\mbox{ times}}{\underbrace{ u, \dots,  u}},
\underset{n-i-1\mbox{ times}}{\underbrace{ v, \dots,  v}},
 u-v)=$$ $$-\int _{s\in\bo} det(\underset{i\mbox{
times}}{\underbrace{ u, \dots,  u}}, \underset{n-i-1\mbox{
times}}{\underbrace{ v, \dots ,v}} ) \frac{\partial}{\partial \nu (s)}
(\alpha\cdot \rho) ds=$$
$$-\int _{s\in\bo} det( u,\dots,  v) \frac{\partial \rho}{\partial
\nu (s)} \cdot \alpha ds.$$ But since $\alpha\leq 0$ and
$\frac{\partial \rho}{\partial \nu (s)}\leq 0$ the last expression
is non-positive. \qed

 {\bf Proof} of Lemma 2.2.3. We have
$$\int_{\Omega} det( u_1,\dots,  u_{n-1}, \beta)=
\int_{\HH^n} \xom det( u_1,\dots,  u_{n-1}, \beta).
$$
By Lemma 2.1.7 the last expression is symmetric with respect to
all the arguments. Hence it is equal to
$$\int _{\HH^n} \beta det ( u_1,\dots,  u_{n-1}, \xom).$$
One easily checks the following
\begin{claim}
$\frac{\partial}{\partial x_i} (\xom)$ is a distribution of order
zero with support on $\bo$. This distribution is equal to
$-\frac{\partial}{\partial x_i}\rfloor vol$.
\end{claim}
Now let us fix a point $s_0\in \bo$. Let us choose an orthonormal
coordinate system $(q_1,\dots, q_n)$ in $\HH^n$, $q_m=t_m+i x_m+j
y_m+k z_m$, such that $\frac{\partial}{\partial t_1}=\nu (s_0)$.

Let $\xi,\, \eta$ be translation invariant vector fields, each of
them parallel to one of the chosen coordinate axes, and at least
one of them is different from $\frac{\partial}{\partial t_1}$. In
the formula for $det ( u_1,\dots,  u_{n-1}, \xom)$
consider the term containing $\xi(\eta(\xom))$. It is a product of
this last term by some smooth function $F$. Let us consider the
integral $\int_{\HH^n} \beta\cdot F\cdot\xi(\eta(\xom))$. We may
assume that at $s_0$ $\eta\in T_{s_0}(\bo)$. Then
$$ \int_{\HH^n}\beta\cdot F\cdot\xi(\eta(\xom))=-\int_{\bo}\beta\cdot
F\cdot\xi(\eta\rfloor vol)= \int_{\bo} \xi(\beta\cdot F)\cdot
(\eta\rfloor vol).$$ But since $\beta|_{\bo} \equiv 0$ the last
expression is equal to $\int _{\bo} \xi(\beta) \cdot F\cdot
(\eta\rfloor vol).$ Note that since $\eta \in T_{s_0}(\bo)$ the
expression under the last integral vanishes at the point $s_0$.
Hence the only summand which remains is
$$\frac {\partial \beta}{\partial \nu (s_0)}\cdot F\cdot(\nu (s_0)\rfloor
vol).$$ It is easy to see that in this case
$$F=det ( u_1|_{\tilde T_{s_0}},\dots,
u_{n-1}|_{\tilde T_{s_0}}),$$ and $ \nu (s_0)\rfloor vol=-ds.$
This proves the lemma. \qed

The next result is a slight generalization of Theorem 2.2.1; it is
completely parallel to Theorem 3.2 of \cite{bedford-taylor}.
\begin{theorem}
Let $\Omega$ be a bounded open set in $\HH^n$. Let $v$ be a
continuous function on $\co$ which is psh in $\Omega$. Let $u$ be
a locally bounded (not necessarily continuous) psh function on
$\Omega$ such that
$$\underset{\zeta \str z \in \bo}{\lim inf} \,(u(\zeta)-v(\zeta))\geq 0;$$
and
$$\underset{\eps \str 0}{\lim} det(\pt u_\eps) \leq det(\pt v)
\mbox{ in } \Omega,$$ where $u_\eps=u\ast\chi_\eps$ and
$\chi_\eps$ is a usual smoothing kernel of psh functions (exactly
as in the complex case, see \cite{lelong}, p.45). Then $u\geq v$
in $\Omega$.
\end{theorem}
{\bf Proof.} Assume that the theorem is false. Then there exists
$z_0 \in \Omega$ such that $u(z_0)<v(z_0)$. Let $\eta
_0=(v(z_0)-u(z_0))/2$. Then for all $0<\eta <\eta_0$ the set
$$G(\eta)=\{z\in \Omega |u(z)+\eta<v(z)\}\ni z_0$$
is nonempty, open (since $u-v$ is upper semi-continuous),
relatively compact subset of $\Omega$ (because of the first
assumption of the theorem).

Let $u_\eps=u\ast\chi_\eps, \,v_\eps=v\ast\chi_\eps$ be
regularizations of $u,v$
so that $u_\eps,
v_\eps$ are defined on
$$ \Omega_{\eps}=\{z\in\Omega|\mbox{ distance from $z$ to $\bo$ exceeds }
\eps\},$$ and $u_\eps \geq u,\, v_\eps \geq v$. Since $v$ is
continuous, $v_\eps \str v$ uniformly on compact subsets of
$\Omega$. Define
$$G(\eta,\delta)=\{z\in\Omega|u(z)+\eta <v(z)+\delta |z-z_0|^2\}.$$
There exists an increasing function $\delta(\eta)>0,\,
0<\eta<\eta_0$, such that $G(\eta,\delta)$ is nonempty, open, and
relatively compact in $\Omega$ for all $0<\delta\leq
\delta(\eta)$. Clearly $z_0\in G(\eta,\delta)$. Next choose
$\eps(\eta,\delta)>0$ so small that $0<\eps<\eps(\eta,\delta)$
implies
$$\Omega_\eps\supset G(\eta/2, \delta), \, 0<\eta<\eta_0,\,
0<\delta<\delta(\eta/2).$$ For such $\eps,\eta, \delta$ let us
define
$$G(\eta,\delta,\eps)=\{z\in G(\eta/2,\delta)|u(z)+\eta< v_\eps(z)
+\delta |z-z_0|^2\}.$$ If $\eps$ is so small that
$|v(z)-v_{\eps}(z)|\leq \eta/4$ whenever $z\in G(\eta/2,\delta)$
and $\eps<\eps(\eta,\delta)$ then it is easy to see that
$$ G(\eta,\delta,\eps)\subset G(3\eta/4,\delta)\subset
G(\eta/2).$$

In particular, $G(\eta,\delta,\eps)$ is a relatively compact
subset of $\Omega_\eps$, so $v_\eps$ is $C^\infty$ in a
neighborhood of the closure of $G(\eta,\delta,\eps)$.

Finally choose $\tau(\eta,\delta,\eps)$ so small that for
$\eta,\delta,\eps$ as above and $0<\tau<\tau(\eta,\delta,\eps)$ we
have that
$$G(\eta,\delta,\eps,\tau):=\{z\in
G(\eta/2,\delta)|u_\tau(z)+\eta<v_\eps(z)+\delta|z-z_0|^2\}$$ is a
nonempty, open, relatively compact subset of $\Omega_\eps$. Since
$u_\tau\geq u$ we have $G(\eta,\delta,\eps,\tau)\subset
G(\eta,\delta,\eps)$, and because $z_0\in G(\eta,\delta,\eps)$ we
have $z_0\in G(\eta,\delta,\eps,\tau)$ for sufficiently small
$\tau$.

We will apply Proposition 2.2.2 with $G(\eta,\delta,\eps,\tau)$
instead of $\Omega$ and the functions defining this set. However
in general this domain does not have smooth boundary. But, by
Sard's lemma, the value $\eta$ is a regular value of the
$C^\infty$- function $v_\eta(z)+\delta|z-z_0|^2- u_\tau(z)$ for
almost all values of $\eta$. Thus we can take sequence of numbers
$\tau_n\str 0$ and apply Proposition 2.2.2  for almost all values
of $\eta$. Consequently we have by Proposition 2.2.2 and Theorem
1.1.17 (iii)
$$\int det(\pt u_\tau)= \int det \pt (u_\tau+\eta) \geq
\int det \pt (v_\eps + \delta |z-z_0|^2)\geq
$$
$$ \int det \pt v_\eps +\delta ^n \int det(\pt |z-z_0|^2)=
\int det \pt v_\eps +\delta^n \cdot c_n
vol(G(\eta,\delta,\eps,\tau)),$$ where all the integrals are taken
over $G(\eta,\delta,\eps,\tau)$, and $c_n$ is a positive constant
depending on $n$ only. When $\tau \str 0$ the open sets
$G(\eta,\delta,\eps,\tau)$ increase to $G(\eta,\delta,\eps)$. If
$\mu=\underset{\tau\str 0}{ \lim} det(\pt u_\tau)$ then we deduce
from the last estimate that
$$\mu (G(\eta,\delta,\eps))\geq \int _{G(\eta,\delta,\eps)}
det(\pt v_\eps) +c_n \delta ^n vol(G(\eta,\delta,\eps))$$ for
almost all $0<\eta<\eta_0,\, 0<\delta<\delta (\eta)$, and $0<\eps
<\eps(\eta, \delta)$. Now let $\eps \str 0$. The measures $det
(\pt v_\eps)$ converge weakly to $\det(\pt v))$ by Theorem 2.1.11.
Also $G(\eta,\delta,\eps)\supset G(\eta,\delta)$. Next
$$\cap_{\eps >0} G(\eta,\delta,\eps)\subset K(\eta,\delta):=
\{z\in \Omega| u(z)+\eta\leq v(z)+\delta|z-z_0|^2\}.$$ Thus for
almost all $\eta$ we have
$$\mu(K(\eta,\delta))\geq \int _{G(\eta,\delta)} det(\pt v) +c_n
\delta^n \cdot vol(G(\eta,\delta)).$$ Let us denote $\nu:=det(\pt
v)$. By assumption $\mu\leq\nu$. Thus we get
$$\nu(K(\eta,\delta))\geq \nu(G(\eta,\delta)) +c_n \delta^n
vol(G(\eta,\delta)).$$  Also $G(\eta,\delta)\subset
K(\eta,\delta)\subset G(\eta ',\delta)$ for $\eta' <\eta$. Hence
$$\nu(G(\eta,\delta))\leq \nu(K(\eta,\delta))\leq \nu(G(\eta ',\delta))
\mbox{ for } \eta' <\eta.$$ However $\eta \mapsto
\nu(G(\eta,\delta))$ is a decreasing function of $\eta$. Hence at
the points of continuity of this function we have
$$\nu(G(\eta,\delta))\geq \nu(G(\eta,\delta))+c_n \delta^n
vol(G(\eta,\delta)).$$ But this contradicts to the fact that
$G(\eta,\delta)$ is a nonempty open set. This proves Theorem 2.2.5
(and hence Theorem 2.2.1) . \qed

\section*{Appendix.}
In this appendix we prove that the linear combinations of
 delta-functions of quaternionic hyperplanes in $\HH^n$
are dense in the space of distributions (this fact was needed in
the proof of Lemma 2.1.7). By the Hahn-Banach theorem it is
equivalent to the injectivity of the Radon transform over
quaternionic hyperplanes. We believe that the injectivity of
quaternionic Radon transform is a well known fact, but we include
the proof for completeness, since we could not find a reference.

Let us fix hyperhermitian metric on $\HH^n$,
i.e. a Euclidean metric such that for any two vectors $x, \, y \in \HH^n$
and any quaternion $a$  with $|a|=1$
$$(x\cdot a, y\cdot a)= (x,\, y).$$
Let $f$ be any smooth compactly supported
function on $\HH^n$.
The quaternionic Radon transform of $f$ is a function on the manifold
of all affine quaternionic hyperplanes defined as
$$Rf(E)= \int_E f(q) dq,$$
where the integration is with respect to the volume form on $E$
defined by the metric.

{\bf Proposition.} {\itshape
The quaternionic Radon transform is injective.
}

{\bf Proof.} We will just present the inversion formula completely
analogous to the complex Radon transform (see \cite{gelfand}). Let
us fix the origin $0\in \HH^n$ for convenience. Let us denote by
${\cal A}$ the manifold of affine quaternionic hyperplanes in
$\HH^n$. For any point $q\in \HH ^n$ let ${\cal P}_q$ denote the
manifold of  quaternionic hyperplanes passing through $q$. For
$E\in {\cal A}$ let us denote by $E^{\perp}$ the quaternionic line
orthogonal to $E$ and passing through the origin 0.

Let us define the  operator
$${\cal D}: C^{\infty}({\cal A}) \str C^{\infty}(\HH ^n) $$
as follows. Let $g\in C^{\infty}({\cal A}) $. Set  $${\cal
D}g(q):= \int_{E\in {\cal P}_q} ( \Delta _{E^\perp}) ^{2(n-1)}
g(E+w)dE,$$ where $\Delta _{E^\perp}$ denotes the (4- dimensional)
Laplacian with respect to $w\in E^\perp$, and the integration is
with respect the Haar measure on ${\cal P}_q$.

{\bf Claim.} {\itshape For any smooth rapidly decreasing function $f$ of $\HH^n$
$${\cal D}(Rf)= c\cdot f,$$
where $c$ is a non-zero constant.
}

It is sufficient to check this claim pointwise, say at $0$. The
operators $R$ and ${\cal D}$ commute with translations and the
action of the group $Sp_n$. Then ${\cal D}(Rf)(0)$ defines a
distribution invariant with respect to the action of $Sp_n$.
Moreover it is easy to check that this distribution is homogeneous
of degree $-4n$ (exactly as the delta-function at $0$). It is easy
to see that there is at most one dimensional space of $Sp_n$-
invariant distributions homogeneous of degree $-4n$. Hence they
must be proportional to the delta-function at $0$. Thus ${\cal
D}(Rf)= c\cdot f$ for some constant $c$. So see that $c\ne 0$ it
is sufficient to check it by an explicit computation for the
function $f(q)=exp(-|q|^2/2)$. \qed

\vskip 0.7cm


\end{document}